\documentclass[12pt]{amsart}
\usepackage{amsthm,amsmath,amssymb,amscd,graphics,psfig}

{\end{list}}
{
   \newtheorem{theorem}{Theorem}[section]
   
   \newtheorem{lemma}[theorem]{Lemma}

   \newtheorem{corollary}[theorem]{Corollary}

}
{\theoremstyle{definition}
   \newtheorem{definition}[theorem]{Definition}

}
\newcommand{\RR}{{\mathbb{R}}}

\newcommand{\CC}{{\mathbb{C}}}

\newcommand{\PP}{{\mathbb{P}}}
\newcommand{\ZZ}{{\mathbb{Z}}}

\newcommand{\cA}{{\mathcal A}}

\newcommand{\cF}{{\mathcal F}}
\newcommand{\cG}{{\mathcal G}}

\newcommand{\cK}{{\mathcal K}}
\newcommand{\cL}{{\mathcal L}}
\newcommand{\cM}{{\mathcal M}}

\newcommand{\res}{\operatorname{res}}
\newcommand{\Id}{\operatorname{Id}}
\newcommand{\Link}{\operatorname{Link}}
\newcommand{\Span}{\operatorname{Span}}

\newcommand{\Star}{\operatorname{Star}}

\newcommand{\Sym}{{\operatorname{Sym}}}

\newcommand{\Tor}{\operatorname{Tor}}

\newcommand{\isomto}{\stackrel{\sim}{\to}}
\newcommand{\isom}{\simeq}

\newcommand{\cb}{\cite{BBFK,BL} }
\newcommand{\LO}{Lefschetz operation }
\newcommand{\HRM}{Hodge-Riemann-Minkowski }
\newcommand{\tensor}{\otimes}
\newcommand{\RRp}{\RR_{\geq 0}}

\begin{document}
\title{Hard Lefschetz Theorem for Nonrational Polytopes}

\author{Kalle Karu}
\thanks{The author was partially supported by NSF grant DMS-0070678}
\address{Department of Mathematics\\ Harvard University\\ 1 Oxford Street\\
Cambridge, MA 02138\\ USA}
\email{kkaru@math.harvard.edu}
% \date{Dec 1, 2001}

\begin{abstract}
The Hard Lefschetz theorem is known to hold for the intersection
cohomology of the toric variety 
associated to a rational convex polytope. One can construct the
intersection cohomology combinatorially 
from the polytope, hence it is well defined even for nonrational
polytopes when there is no variety associated to it. We
prove the Hard Lefschetz theorem for the intersection cohomology of a
general polytope. 
\end{abstract}

\maketitle

\addtocounter{section}{-1}

\section{Introduction}
If $P$ is a simple convex $n$-dimensional polytope in $\RR^n$, define 
\[ h_k = \sum_{i\geq k} f_i (-1)^{i-k}\binom{i}{k},\]
where $f_i$ is the number of $i$-dimensional faces of $P$. Then
knowing the numbers $h_k$ is equivalent to knowing the numbers $f_i$.
The following conditions are satisfied by the numbers $h_k$:
\begin{enumerate}
\item Dehn-Sommerville equations: $h_i = h_{n-i}$.
\item Unimodality: $1= h_0\leq h_1\leq \ldots\leq h_{[n/2]}$.
\end{enumerate}
To prove the unimodality condition (in fact a stronger condition
conjectured by P.~McMullen), R.~Stanley in \cite{S1}
constructed a projective toric variety $X_P$ so that $h_i$
are the even Betti numbers of the (singular) cohomology of $X_P$. The two
conditions then follow from the Poincar\'e duality and the Hard
Lefschetz theorem for $X_P$. Later in \cite{S2} Stanley generalized
the definition of the $h$-vector to an arbitrary polytope $P$ in such
a way that if the polytope is rational the numbers $h_i$ are the
intersection cohomology Betti numbers of the associated toric
variety (see \cite{F} for a validation of this claim). Stanley also
proved that the Dehn-Sommerville equations hold for the 
generalized vector $h$. If the polytope $P$ is rational, the
unimodality condition follows from the Hard Lefschetz theorem in
intersection cohomology of the associated toric variety. The goal of
this paper is to prove the Hard Lefschetz theorem for a general
polytope, generalizing McMullen's purely combinatorial proof of the
Hard Lefschetz theorem for a simple polytope \cite{M}. 

Our proof is based on the description of the intersection cohomology
of a fan (e.g., the normal fan of a polytope) given in
\cb. When the fan is simplicial, the torus-equivariant cohomology
of the corresponding toric variety can be identified with the space of
conewise
polynomial functions on the fan. This set is naturally a module over
the algebra $A$ of global polynomial functions. Taking the quotient by
the ideal $I\subset A$ generated by global linear functions gives us
the non-equivariant intersection cohomology.
 In order to generalize this, one views the
conewise polynomial functions as global sections of a sheaf on the
fan with respect to a certain topology. On a general fan $\Sigma$ one
defines the equivariant intersection cohomology sheaf $\cL_\Sigma$; the
equivariant 
intersection cohomology of the fan $\Sigma$ is then the $A$-module of
global sections of this sheaf. Taking the quotient by the ideal $I$,
we obtain the non-equivariant
intersection cohomology $IH(\Sigma)$ of the fan $\Sigma$.

If the fan $\Sigma$ is projective, there exists a conewise linear
strictly convex function $l$ on $\Sigma$. Multiplication with this
function defines a degree $2$ homomorphism of the sheaf $\cL$ (we
assume that linear functions have degree 2), and it also induces a
degree $2$ map of $IH(\Sigma)$. 

\begin{theorem}[Hard Lefschetz] Let $\Sigma$ be a projective fan of
dimension $n$ and let $l$ be a conewise linear strictly convex
function on $\Sigma$. 
Then multiplication with $l$ defines a \LO on $IH(\Sigma)$. That means
\[ l^i: IH^{n-i}(\Sigma)\longrightarrow IH^{n+i}(\Sigma) \]
is an isomorphism for all $i>0$.
\end{theorem}

The theorem is known for simplicial fans \cite{S1, M} and
for rational fans \cite{S2} (a fan in $\RR^n$ is rational if 
every cone is generated by a set of vectors with rational
coordinates). V. Timorin \cite{T1} recently proved the Hard Lefschetz
theorem for 
certain fans in which all cones of dimension less than $n$ are
simplicial. 

The Hard Lefschetz theorem follows from an even stronger
statement, the 
Hodge-Riemann-Minkowski relations. Using the bilinear Poincar\'e
duality pairing \cite{BBFK}:
\begin{alignat*}{2}
 IH^{n-i}(\Sigma) & \times IH^{n+i}(\Sigma) & \longrightarrow  & \RR \\
 (x & , y) &\longmapsto & \langle x  \cdot y\rangle,
\end{alignat*}
we define the quadratic form $Q_l(x)=\langle l^{i}x\cdot x\rangle$ on
$IH^{n-i}(\Sigma)$. Then if 
\[ IP^{n-i}(\Sigma) = \ker l^{i+1}: IH^{n-i}(\Sigma) \longrightarrow
IH^{n+i+2}(\Sigma)\]
is the primitive cohomology, we have

\begin{theorem}[Hodge-Riemann-Minkowski bilinear relations]\label{thm-HRM}
The quadratic form 
\[ (-1)^{\frac{n-i}{2}} Q_l \]
is positive definite on $IP^{n-i}(\Sigma)$ for all $i\geq 0$.
\end{theorem}

We remark that the intersection cohomology of a complete fan is
nonzero only in even degrees, hence the power ${\frac{n-i}{2}}$ is
always an integer. McMullen in \cite{M} proved the Hodge-Riemann-Minkowski
relations for simplicial fans (see also the proof by Timorin \cite{T2}). The
names Hodge-Riemann refer to the usual Hodge-Riemann bilinear
relations for K\"ahler manifolds; the inequalities in degree $n-i = 2$
reduce to the Minkowski inequalities on polytopes \cite{M}.

We will assume the Hodge-Riemann-Minkowski inequalities in
the simplicial case and deduce Theorem~\ref{thm-HRM} from it. 

{\bf Acknowledgments.} I would like to thank Gottfried Barthel, Valery
Lunts and G\"unter Ziegler for helpful comments. The preprint in the current
form owes much to the advice of Karl-Heinz Fieseler. Almost the entire
last section has been rewritten following his suggestions. I take
credit for any remaining errors.

\section{Intersection cohomology of fans}

We start by recalling the definition of the intersection cohomology of a
fan. The main references here are \cite{BBFK,BL,B}.  

Let $V$ be a real vector space of dimension
$n$. We denote by $A$ the
algebra of polynomial functions on $V$:
\[ A = \Sym V^*,\]
graded so that linear functions have degree $2$. The ideal generated
by linear functions in $A$ is denoted by $I$, and the quotient of an
$A$-module $M$ by the submodule $I M$ is denoted by $\overline{M}$.

\subsection{Cones and fans} 
 A {\em polyhedral cone} $\sigma$ in $V$ is a cone generated by a
finite set of vectors $v_1,\ldots v_N \in V$:
\[ \sigma = \RR_{\geq 0} v_1+\ldots +\RR_{\geq 0} v_N.\]
We only consider pointed cones, which  means that $\{0\}$ is the
largest sub-vector space in a cone. 
Given a cone
$\sigma\subset V$ and a linear function $f$ on $V$ such that $f$ takes
non-negative 
values on $\sigma$, the cone $\tau=\sigma\cap \ker{f}$ is called a face of
$\sigma$. We denote this relation by $\tau\prec\sigma$. A face of
dimension $\dim\sigma-1$ is called a facet of $\sigma$.

A {\em fan} $\Sigma$ in $V$ is a finite collection of cones in $V$
such that (i) if $\sigma\in\Sigma$ and $\tau$ is a face of $\sigma$
then $\tau\in\Sigma$; and (ii) if $\sigma_1, \sigma_2\in\Sigma$ then
the intersection $\sigma_1\cap\sigma_2$ is a face of both $\sigma_1$
and $\sigma_2$. 
Examples of fans include the fan $[\sigma]$ consisting of all faces of
a cone $\sigma$, and the fan $\partial\sigma =
[\sigma]-\{\sigma\}$. Given a cone $\sigma\in\Sigma$ we define
\begin{align*}
\Star_\Sigma(\sigma) &= \{\delta\in\Sigma\mid\sigma\prec\delta\}, \\
\overline{\Star}_\Sigma(\sigma) &= \{\tau\in\Sigma\mid\tau\prec\delta
\text{ for some } \delta\in \Star_\Sigma(\sigma)\}, \\
\Link_\Sigma(\sigma) &= \{\nu\in\overline{\Star}_\Sigma(\sigma) \mid
\nu\cap\sigma = \{0\}\}.
\end{align*}
Then $\overline{\Star}_\Sigma(\sigma)$ and $\Link_\Sigma(\sigma)$ are
subfans of $\Sigma$.

The {\em support} of a fan $\Sigma$ is $|\Sigma| =
\cup_{\sigma\in\Sigma} \sigma$. A fan is {\em complete} if its support
is $V$. The {\em relative interior} of a cone $\sigma$ is
$\sigma- |\partial\sigma|$.

A cone $\sigma$ is called {\em simplicial} if it can be generated by
$\dim \sigma$ vectors. A fan is simplicial if all its cones are simplicial.

\subsection{Sheaves on a fan}

We consider a fan $\Sigma$ as a set of cones and  define a topology on
it so that open sets are the subfans of $\Sigma$. We can then consider
sheaves of $\RR$-vector spaces on 
$\Sigma$ with respect to this topology. To give a sheaf $\cF$ on
$\Sigma$, it suffices to specify for each cone $\sigma\in\Sigma$ the stalk
$\cF_\sigma = \cF([\sigma])$ and the restriction map
$\cF_\sigma \to \cF(\partial\sigma)$. 

The {\em structure sheaf} $\cA_\Sigma$ is defined by setting
$\cA_{\Sigma,\sigma}= \Sym (\Span{\sigma})^*$, the space of polynomial
functions on $\sigma$; the restriction map
$\cA_{\Sigma,\sigma}\to \cA_{\Sigma}(\partial\sigma)$ is defined by
restriction of 
functions. Then $\cA_\Sigma$ is a sheaf of algebras, naturally graded
by degree  (again, conewise linear functions have degree
$2$). Multiplication with elements of $A$ gives  $\cA_\Sigma$ the
structure of a sheaf of $A$-modules. 

An {\em equivariant intersection cohomology sheaf} $\cL_\Sigma$ of
$\Sigma$ is a sheaf of $\cA_\Sigma$-modules satisfying the following
properties:
\begin{enumerate}
\item Normalization: $\cL_{\Sigma,0} = \RR$.
\item Local freeness: $\cL_{\Sigma,\sigma}$ is a free
$\cA_{\Sigma,\sigma}$-module for any $\sigma\in\Sigma$.
\item Minimal flabbiness: Modulo the ideal $I\subset A$ the restriction map
induces an isomorphism
\[ \overline{\cL}_{\Sigma,\sigma} \longrightarrow
\overline{\cL_\Sigma(\partial\sigma)}.\] 
\end{enumerate}
Note that the third property of an equivariant intersection
cohomology sheaf implies that $\cL_\Sigma$ is flabby: the
restriction map $\cL_\Sigma(U)\to \cL_\Sigma(V)$ is surjective for any
open sets $V\subset U\subset \Sigma$. Equivariant intersection
cohomology sheaves exist for any fan $\Sigma$, and any two of them
are isomorphic, hence we may call a sheaf of  $\cA_\Sigma$-modules
satisfying the three properties {\em the} equivariant intersection
cohomology sheaf. As a special case, $\cL_\Sigma\isom \cA_\Sigma$ if
and only if the fan $\Sigma$ is simplicial.

One can 
construct the sheaf $\cL_\Sigma$ by induction on the dimension of cones
as follows: assume that $\cL_\Sigma$ is defined on $\partial\sigma$,
and put 
\[ \cL_{\Sigma,\sigma} =
\overline{\cL_\Sigma(\partial\sigma)}\tensor \cA_{\Sigma,\sigma}\]
with the obvious restriction maps. (Note! All tensor products are
taken over $\RR$ unless specified otherwise.)
From this construction it is clear that if $\Delta\subset \Sigma$ is a
subfan, then $\cL_\Delta \isom \cL_\Sigma|_\Delta$. To simplify notation,
we will often omit the subscript $\Sigma$ and call the equivariant
intersection cohomology sheaf simply
$\cL$. For example, the stalk $\cL_\sigma$ for a cone $\sigma$ is
defined up to 
isomorphism independent of the fan containing $\sigma$. We use similar
omission of the subscript for the structure sheaf $\cA_\Sigma$. Then
$\cA_\sigma = \cA_{\Sigma,\sigma}$ is the algebra of polynomial
functions on $\sigma$. 

The (non-equivariant) {\em intersection cohomology} of a fan $\Sigma$
is defined to be the $A$-module of global sections of the equivariant
intersection cohomology sheaf modulo the ideal $I$:
\[ IH(\Sigma) = \overline{\cL(\Sigma)}\]
Since the equivariant intersection cohomology sheaf is a graded
$A$-module and $I$ is a graded ideal, $IH(\Sigma)$ inherits a natural
grading.

\subsection{Morphisms of fans} Let $\Sigma$ be a fan in $V$ and
$\Delta$ a fan in $W$. A {\em conewise linear morphism}
$\Sigma\to\Delta$ is a pair 
$(|\phi|,\phi)$, where  $|\phi|: |\Sigma|\to |\Delta|$  maps each cone
$\sigma\in\Sigma$ linearly to some cone $\delta\in\Delta$, and
$\phi$ is the induced continuous map between the fan spaces. We often denote
the pair $(|\phi|,\phi)$ simply by $\phi$. Note that we do not require
the map $|\phi|$ to be the restriction of a linear map $V\to W$;
rather, the map $|\phi|$ is linear only when restricted to a cone.

If $\phi:\Sigma\to\Delta$ is a conewise linear morphism then by composing a
section of $\cA_\Delta$ with $|\phi|$, we get a section of
$\cA_\Sigma$, hence $\phi$ defines a morphism of ringed spaces. If
$\cG$ is a sheaf of $\cA_\Sigma$-modules then 
$\phi_*\cG$ has the structure of an $\cA_\Delta$-module. For a
sheaf $\cF$ on 
$\Delta$ we denote the pullback sheaf by $\phi^{-1}\cF$. If, moreover,
$\cF$ is a sheaf of $\cA_\Delta$-modules, we define the sheaf of
$\cA_\Sigma$-modules  
\[ \phi^*\cF = \phi^{-1}\cF \otimes_{\phi^{-1}\cA_\Delta}
\cA_\Sigma.\]

As an example, suppose $\phi:\Sigma\to\Delta$ is a conewise linear
isomorphism. Then 
since $|\phi|$ maps every cone $\sigma\in\Sigma$ isomorphically onto
some cone $\delta\in\Delta$, it follows that $\phi_* \cL_\Sigma\isom
\cL_\Delta$. In particular, we get an isomorphism of the vector spaces of
global sections $\cL(\Sigma)\isom \cL(\Delta)$; however, the
$A$-module structures may be very different and the intersection
cohomologies of the two fans need not be isomorphic. If we know that
the spaces of global sections are finitely generated free $A$-modules
(for example, when $\Sigma$ and $\Delta$ are complete), then by counting
the dimensions of the graded pieces, we get a non-canonical
isomorphism between the intersection cohomologies. 

\subsection{Subdivisions of fans} A conewise linear morphism of fans
$s:\Sigma\to\Delta$ in the same space $V$ is called a {\em
subdivision} if $|\Sigma| = |\Delta|$ and $|s|$ is the identity map.  
It follows from the decomposition theorem \cite{BBFK,BL} that if $s$ is a subdivision
then $\cL_\Delta$ is a direct summand of $s_*\cL_\Sigma$ as sheaves of
$\cA_\Delta$-modules; in particular, $IH(\Sigma)$ contains
$IH(\Delta)$ as a direct summand. 

If $\Delta$ is an arbitrary fan, there exists a subdivision
$s:\Sigma\to\Delta$ such that $\Sigma$ is simplicial, for example, the
barycentric subdivision constructed below. Since
$\cL_\Sigma\isom \cA_\Sigma$, using the decomposition theorem we can
consider sections of $\cL_\Delta$ as conewise polynomial functions on
the subdivision $\Sigma$. 

Now let us construct a simplicial subdivision of an arbitrary fan
$\Sigma$. First, suppose $\sigma\in\Sigma$ is a cone of dimension at
least $2$ such that any
other cone $\tau\in\Sigma$ containing $\sigma$ can be written as a
direct sum $\tau=\sigma\oplus \rho$ for some $\rho\in
\Link_\Sigma(\sigma)$. If $\RR_{\geq 0} v$ is a ray intersecting the
relative interior of $\sigma$, we define the {\em star subdivision} of
$\Sigma$ at $\RR_{\geq 0} v$ to be the fan
\[ (\Sigma - \Star_\Sigma(\sigma)) \cup \{\RR_{\geq 0} v + \tau + \rho
\mid \tau\in\partial\sigma, 
\rho\in\Link_\Sigma(\sigma)\}.\]
To simplify notation, we assume that the star subdivision of $\Sigma$
at $\RRp v$ where $v$ lies in a cone of dimension less than $2$ does
not change the fan.

Given an arbitrary fan $\Sigma$, for every cone $\sigma\in\Sigma$ let
$v_\sigma$ be a vector in the relative interior of $\sigma$.
We can perform a sequence of star 
subdivisions at $\RR_{\geq 0} v_\sigma$, starting with cones $\sigma$ of
maximal dimension, then cones of one smaller dimension, and so on. The
resulting subdivision is simplicial and it is called a {\em
barycentric subdivision} of $\Sigma$.

\subsection{Evaluation map in the top degree cohomology}

Following M.~Brion \cite{B} we construct an isomorphism $IH^{2n}(\Sigma) \to
\RR$ for a complete $n$-dimensional fan $\Sigma$. 

First assume that $\Sigma$ is a simplicial fan. For every
$n$-dimensional cone $\sigma\in\Sigma$ we choose linear functions
$H_{\sigma,1},\ldots,H_{\sigma,n}$ defining the facets of $\sigma$,
such that $H_{\sigma,i}$ are non-negative on $\sigma$ and their wedge
product has length 1 with respect to some metric on $\Lambda^n V^*$. Let
$\Phi_\sigma$ be the product of the $H_{\sigma,i}$. Thus, up to a
constant factor, $\Phi_\sigma$ is the unique degree $n$ polynomial that
vanishes on $\partial \sigma$; the constant is determined by the
metric on $\Lambda^n V^*$ together with the requirement that
$\Phi_\sigma$ is non-negative on $\sigma$. 

We identify $\cL_\Sigma = \cA_\Sigma$ and consider a section
$f\in\cL(\Sigma)$ as a conewise polynomial function. Then for a
maximal cone $\sigma$, the stalk $f_\sigma \in \cA_\sigma = A$ is
a polynomial and we can define a rational function 
\[ \langle f\rangle = \sum_{\sigma\in\Sigma,\; \dim\sigma = n}
\frac{f_\sigma}{\Phi_\sigma}.\]
Brion proved that the poles of
this rational function cancel out so that the result is in
fact a polynomial, thus defining a degree $-2n$ map of $A$-modules
\[ \langle \cdot\rangle: \cL(\Sigma) \longrightarrow A. \]
Since $\langle \cdot\rangle$ maps $\cL(\Sigma)^{2n-2}$ to zero, it
induces a map
\[ \langle \cdot\rangle: IH^{2n}(\Sigma) \longrightarrow A^0 = \RR.\]
Because this map is surjective (by considering a function $f$ supported
on $\sigma$ such that $f|_\sigma = \Phi_\sigma$) and $\dim
IH^{2n}(\Sigma) = \dim IH^{0}(\Sigma) = 1$ (see the next section about
Poincar\'e duality), it follows that the map $\langle \cdot\rangle$ is an
isomorphism. 

Now suppose that $\Delta$ is an arbitrary complete fan. We choose a
simplicial subdivision $s: \Sigma\to\Delta$ and, using the
decomposition theorem, embed
\[ \cL_\Delta \subset s_* \cL_\Sigma.\]
Restricting the map $\langle \cdot\rangle$ defined on the global sections of $\cL_\Sigma$ to
the global sections of  $\cL_\Delta$ gives us a map in cohomology
\[ \langle \cdot\rangle: IH^{2n}(\Delta) \hookrightarrow IH^{2n}(\Sigma) \longrightarrow \RR.\]
Since both top degree cohomology spaces have dimension 1, this
composition is an isomorphism. 

We sometimes use a subscript to
indicate the fan for which the evaluation map is constructed: 
\[ \langle \cdot\rangle_\Delta: IH^{2n}(\Sigma)\longrightarrow\RR.\]

\subsection{Fans with boundary and the Poincar\'e pairing}

As in \cite{BBFK} we consider a (possibly noncomplete) fan $\Sigma$
such that all maximal cones of $\Sigma$ have dimension $n$. Then the
boundary $\partial\Sigma$ of $\Sigma$ consists of all faces of those
$n-1$ dimensional cones in $\Sigma$ that are contained in only one $n$
dimensional cone. We further restrict ourselves to the case of {\em
  quasi-convex} fans \cite{BBFK} where the
boundary $\partial \Sigma$ is a real homology manifold.
The examples we wish to consider are complete fans (with
$\partial\Sigma=\emptyset$), and fans of the type $\Sigma =
\overline{\Star}_\Delta(\delta)$ 
for some cone $\delta$ in a complete fan $\Delta$.

If $\Sigma$ is a quasi-convex fan with boundary $\partial\Sigma$, let
$\cL(\Sigma,\partial\Sigma)\subset \cL(\Sigma)$ be the space of
sections vanishing on the boundary, and let 
\[ IH(\Sigma,\partial\Sigma) =
\overline{\cL(\Sigma,\partial\Sigma)}.\]
(Note: in \cite{BBFK} the notation for $IH(\Sigma)$ and
$IH(\Sigma,\partial\Sigma)$ is reversed compared to the one used here.)

It is proved in \cite{BBFK} that for a quasi-convex fan $\Sigma$ both
$\cL(\Sigma)$ and $\cL(\Sigma,\partial\Sigma)$ are finitely generated
free $A$-modules. Moreover, there exists a bilinear
nondegenerate map, called Poincar\'e pairing,
\[ IH^{n-k}(\Sigma) \times IH^{n+k}(\Sigma,\partial\Sigma) \longrightarrow\RR,
\qquad k\in\ZZ.\]
The pairing is constructed by considering sections of $\cL_\Sigma$ as
conewise polynomial functions on a simplicial subdivision $\Delta$ of
$\Sigma$, multiplying these functions, and composing with the
evaluation map $\langle \cdot\rangle:
IH^{2n}(\Delta) \longrightarrow \RR$. We denote this pairing by
\[ (h,g) \longmapsto \langle h\cdot g\rangle.\] 

It is clear that the Poincar\'e pairing depends on the choice of a
subdivision $\Delta$ and the embedding
$\cL(\Sigma)\subset\cL(\Delta)$. We will fix one set of these choices in
Section~4 and call the corresponding  pairing {\em the} Poincar\'e pairing.

Another construction we will use later is the following. Let $\Sigma$
be a complete fan and $\Sigma_1,\Sigma_2$ subfans of $\Sigma$ such
that 
\[ \Sigma_1\cup\Sigma_2 = \Sigma,\qquad \Sigma_1\cap\Sigma_2 =
\partial\Sigma_1=\partial\Sigma_2. \]
Assume that $\Sigma_1$, hence also $\Sigma_2$, is quasi-convex. Then
we have an exact sequence of free $A$ modules
\[
0\longrightarrow\cL(\Sigma_2,\partial\Sigma_2)\longrightarrow\cL(\Sigma)\longrightarrow\cL(\Sigma_1)\longrightarrow 0,
\]
which splits non-canonically. Taking quotients by the ideal $I$ we get
an exact sequence of cohomology spaces 
\[ 0\longrightarrow IH(\Sigma_2,\partial\Sigma_2)\longrightarrow IH(\Sigma)\longrightarrow IH(\Sigma_1)\longrightarrow
0.\]
A splitting of the first sequence gives us a splitting of the second one.

\section{Projective Fans}

Let $\Sigma$ be a complete $n$-dimensional fan and
$l\in\cA^2(\Sigma)$ a conewise 
linear function on $\Sigma$. For an $n$-dimensional
cone $\sigma\in\Sigma$ we consider $l_\sigma = l|_\sigma$ as an
element in $A^2$. 

\begin{definition}\label{def-proj}  A function $l\in\cA^2(\Sigma)$ is
called {\em strictly 
convex} if $l_\sigma(v)<l(v)$ for any maximal cone $\sigma\in\Sigma$
and any  $v\notin\sigma$. A complete fan $\Sigma$ is called {\em
projective} if there exists a strictly convex function $l\in\cA^2(\Sigma)$.
\end{definition}

We note that if $f\in A$ is any linear function then $l$ is strictly
convex if and only if $l+f$ is strictly convex. If $l$ is strictly
convex, then modifying $l$ by some $f\in A$ if necessary, we may
assume that $l(v)>0$ for $v\neq 0$. Then the set $\{v\mid l(v)\leq
1\}$ is a convex $n$-dimensional polytope, with faces corresponding to
cones in $\Sigma$.  

The main goal of this paper is to prove that if $l\in\cA^2(\Sigma)$
is a strictly convex function on a complete fan then multiplication
with $l$ induces an isomorphism 
\[ l^k: IH^{n-k}(\Sigma)\longrightarrow IH^{n+k}(\Sigma).\]

\begin{definition} Let $H$ be a finite dimensional graded vector space
and $l:H\to H$ a degree $2$ linear map. We say that $l$ is a {\em \LO
centered at degree $n$} if  
\[ l^k: H^{n-k}\longrightarrow  H^{n+k}\]
is an isomorphism for all $k>0$.
\end{definition}

If $l:H\to H$ is a \LO centered at degree $n$, we define the primitive
part of $H$ with respect to $l$:
\[ P_l^{n-k} = H^{n-k}\cap\ker l^{k+1}.\]
Then if $\{p_i\}_{i=1}^N$ is a homogeneous basis of $P_l$, one can construct a
basis of $H$ of the form $\{ l^{j_i} p_i\mid1\leq i\leq N, 0\leq j_i\leq
n-\deg(p_i)\}$. 

% Another way to view a \LO is to extend it to an action
% of the Lie group $sl_2$ on the (complexified) vector space $H$ (see
% \cite{GH}). 

\subsection{Hodge-Riemann-Minkowski bilinear relations}

Let $\Sigma$ be a complete $n$-dimensional fan and $l\in\cA^2(\Sigma)$
a strictly 
convex function on $\Sigma$. Using the Poincar\'e pairing we define a
quadratic form $Q_l$ on $IH^{n-k}(\Sigma)$:
\[ Q_l(h) = \langle l^k h \cdot h\rangle. \]

\begin{definition} We say that $Q_l$ satisfies
{\em Hodge-Riemann-Minkowski bilinear relations} if
for any $k\geq 0$ the quadratic form 
\[ (-1)^\frac{n-k}{2} Q_l \]
is positive definite when restricted to the primitive cohomology 
\[ IP_l^{n-k}(\Sigma) = IH^{n-k}(\Sigma) \cap\ker l^{k+1}.\] 
\end{definition}

It is clear that if $Q_l$ satisfies the Hodge-Riemann-Minkowski bilinear
relations then $l^k: IH^{n-k}(\Sigma)\to IH^{n+k}(\Sigma)$ is
injective for all $k>0$. Poincar\'e duality then implies that these
maps are isomorphisms and $l$ defines a \LO on 
$IH(\Sigma)$ centered at degree $n$. 

It is also useful to write the
Hodge-Riemann-Minkowski bilinear relations in terms of the signature
of $Q_l$.

\begin{lemma} $Q_l$ satisfies Hodge-Riemann-Minkowski bilinear
relations if and only if for any $k\leq n$ the signature of the
quadratic form $Q_l$ on $IH^k(\Sigma)$ is 
\[ (\sum_{0\leq j\leq k,\; 4|j} h^j-h^{j-2}, \sum_{0\leq j\leq
k,\; 4 \nmid j} h^j-h^{j-2}),\]
where  $h^j = \dim IH^j(\Sigma)$.
\end{lemma}

{\bf Proof.} Note that each of the two statements whose equivalence we are
proving implies that $Q_l$ is nondegenerate on
$IH^k(\Sigma)$. Hence we may assume that $l$ defines a \LO and, in
particular, that $l:IH^{k-2}(\Sigma) \to
IH^{k}(\Sigma)$ is injective for $k\leq n$. Thus we have a decomposition
\[ IH^{k}(\Sigma) = l\cdot IH^{k-2}(\Sigma) \oplus IP_l^{k} (\Sigma). \]

Since both the Poincar\'e pairing and the operation $l$
are defined by multiplication of functions, it follows that $l$ is
self-adjoint with respect to the pairing:
\[ \langle lx\cdot y\rangle = \langle x\cdot ly\rangle.\]
This implies that for $k\leq n$ we have
\begin{enumerate}
\item[(a)] $Q_l(lx) = Q_l(x)$ for any $x\in IH^{k-2}(\Sigma)$.
\item[(b)] The primitive cohomology $IP_l^{k}(\Sigma)$ is 
orthogonal to $l \cdot IH^{k-2}(\Sigma)$ with respect to the
bilinear form $B_l$ on $IH^{k}(\Sigma)$:
\[ B_l(x,y) = \langle l^k x\cdot y\rangle.\]
\end{enumerate}
Thus the decomposition $IH^{k}(\Sigma) = l\cdot IH^{k-2}(\Sigma)
\oplus IP_l^{k} (\Sigma)$ is orthogonal with respect to the bilinear
form $B_l$ and the signature of $Q_l$ on $IH^{k}(\Sigma)$ is the sum of
the signatures of $Q_l$ on $IH^{k-2}(\Sigma)$ and on $IP_l^{k}
(\Sigma)$. The claim of the lemma now follows by induction on $k$.
\qed

McMullen in \cite{M} proved that the \HRM relations hold if the fan is
simplicial (or equivalently, the corresponding polytope is simple). His proof,
however, was written in terms of the polytope algebra and the weight
spaces of a polytope. In the following, we want to assume
that the \HRM relations hold 
on the intersection cohomology of a simplicial fan when the quadratic
form is defined by multiplication of functions as above. 
There are two ways we can make this assumption valid. First, if the
fan is rational, then the algebra $\cA(\Sigma)$ is isomorphic to the
equivariant cohomology $H_T(X_\Sigma)$ of the corresponding toric
variety (e.g., see \cite{GKM}). Then the \HRM relations follow from
the corresponding relations on $H_T(X_\Sigma)$. The case of a
nonrational simplicial fan follows by a small deformation to a
rational fan, as in the original proof of Stanley \cite{S1}.
Second, if we wish to stay in the realm of algebra and combinatorics,
we can apply the Koszul functor \cite{GKM} to the cellular complex
\cite{BBFK, BL}
computing the global sections $\cA(\Sigma)$. As a result we recover 
the weight spaces together with the product rule used by
McMullen. Since this calculation is straight-forward, we leave it to
the reader.

\section{Products of fans}

\subsection{The K\"unneth theorem}

We consider fans $\Sigma$ in $V$ and $\Delta$ in $W$, and their
product fan $\Sigma\times\Delta$ in $V\times W$. Let $\pi_1:V\times
W\to V$ and $\pi_2:V\times W\to W$ be the two projections, inducing
the maps of fans $\pi_1:\Sigma\times\Delta\to\Sigma$ and
$\pi_2:\Sigma\times\Delta\to\Delta$. Then it is 
easy to see that 
\[ \pi_1^{-1}\cA_\Sigma \tensor \pi_2^{-1}\cA_\Delta \isom
\cA_{\Sigma\times\Delta}.\] 
Indeed, there is a natural morphism $\pi_1^{-1}\cA_\Sigma \tensor
\pi_2^{-1}\cA_\Delta \to \cA_{\Sigma\times\Delta}$ defined by pull-back
of functions. On each cone $\sigma\times\delta\in\Sigma\times\Delta$ 
this map is an isomorphism  $\cA_\sigma\tensor \cA_\delta \isomto
\cA_{\sigma\times\delta}$, hence it is an isomorphism of sheaves. Using 
this isomorphism, we can view $\pi_1^{-1}\cL_\Sigma \tensor
\pi_2^{-1}\cL_\Delta$ as a sheaf of $\cA_{\Sigma\times\Delta}$ modules.

\begin{theorem}\label{thm-kunneth} Let $\Sigma$ be a fan in $V$ and
$\Delta$ a fan in 
$W$. Then 
\[ \pi_1^{-1}\cL_\Sigma \tensor \pi_2^{-1}\cL_\Delta \isom
\cL_{\Sigma\times\Delta}.\]
Moreover, there exists an isomorphism of $\Sym (V\times W)^* \isom \Sym
V^*\otimes \Sym W^*$ modules
\[ \cL(\Sigma) \tensor \cL(\Delta) \isom
\cL(\Sigma\times\Delta).\]
\end{theorem}

{\bf Proof.} We prove both statements by induction on the dimension of
the fan $\Sigma\times\Delta$, the case where either $\Sigma$ or
$\Delta$ has dimension $0$ being trivial. Let us denote $A=\Sym V^*$
and $B=\Sym W^*$.

We start by proving the second statement, assuming the first. We have 
\[ (\pi_1)_* \cL_{\Sigma\times\Delta} \isom (\pi_1)_*
(\pi_1^{-1}\cL_\Sigma \tensor \pi_2^{-1}\cL_\Delta) \isom \cL_\Sigma
\otimes \cL(\Delta).\]
Taking global sections of this push-forward gives the second
statement. 

Next we prove the first statement of the theorem, assuming the theorem
for fans of smaller dimension. We have to show that the sheaf $\cM =
\pi_1^{-1}\cL_\Sigma \tensor \pi_2^{-1}\cL_\Delta$ satisfies the three
conditions in the definition of an equivariant intersection cohomology
sheaf. Note that for any cone $\sigma\times\delta
\in\Sigma\times\Delta$ we have 
\[ \cM_{\sigma\times\delta} = \cL_\sigma\otimes\cL_\delta.\]
From this we deduce the normalization condition $\cM_{0\times 0} =
\RR\otimes\RR = \RR$ and the local freeness of $\cM$: since
$\cL_\sigma, \cL_\delta$ are free $\cA_\sigma, \cA_\delta$ modules,
respectively, it follows 
that $\cL_\sigma\otimes\cL_\delta$ is a free
$\cA_\sigma\otimes\cA_\delta$ module.

It remains to show that the restriction map $\cM_{\sigma\times\delta}\to
\cM(\partial(\sigma\times\delta))$ induces an isomorphism after taking
quotients by the ideal $I$. From the formula 
\[ \partial(\sigma\times\delta) = \partial\sigma \times \delta
\bigcup_{\partial\sigma\times\partial\delta}
\sigma\times\partial\delta\]
we get that $\cM(\partial(\sigma\times\delta))$ is the kernel of the
map 
\[  \cM(\partial\sigma \times \delta) \oplus
\cM(\sigma\times\partial\delta)\stackrel{\bigl(\begin{smallmatrix} \res\\
    -\res\end{smallmatrix}\bigr)}{\longrightarrow}
\cM(\partial\sigma\times\partial\delta),\]
defined by restrictions with appropriate signs. Also note that the
restriction map $\cM_{\sigma\times \delta}\to
\cM(\partial(\sigma\times\delta))$ is defined by 
\[ \cM_{\sigma\times\delta} \stackrel{\bigl(\begin{smallmatrix} \res\\
    \res\end{smallmatrix}\bigr)}{\longrightarrow} 
\cM(\partial\sigma \times \delta) \oplus
\cM(\sigma\times\partial\delta).\]
We now apply induction assumption to the fans $\partial\sigma \times
\delta$, $\sigma\times\partial\delta$ and
$\partial\sigma\times\partial\delta$ to get an exact sequence
\[ 0\longrightarrow \cM(\partial(\sigma\times\delta)) \longrightarrow \cL(\partial\sigma)\tensor
\cL_\delta \oplus \cL_\sigma \tensor \cL(\partial\delta)\longrightarrow
\cL(\partial\sigma)\tensor\cL(\partial\delta)\longrightarrow 0,\]
where exactness on the right follows from the flabbiness of $\cL$. The
restriction map from $\cM_{\sigma\times \delta} =
\cL_\sigma\otimes\cL_\delta$ to $\cM(\partial(\sigma\times\delta))$ is
defined by a map to the middle term of this sequence.
Assuming for a moment that this sequence remains exact after tensoring
with the $A\tensor B$-module $\RR$, and using that the restriction
maps induce isomorphisms
$\overline{\cL}_\sigma\isom \overline{\cL(\partial\sigma)}$ and
$\overline{\cL}_\delta\isom \overline{\cL(\partial\delta)}$, we get an
exact sequence
\[ 0\to \overline{\cM(\partial(\sigma\times\delta))} \longrightarrow
\overline{\cL(\partial\sigma)}\tensor 
\overline{\cL(\partial\delta)} \oplus \overline{\cL(\partial\sigma)}
\tensor \overline{\cL(\partial\delta)}\stackrel{\bigl( \begin{smallmatrix} \Id\\
    -\Id\end{smallmatrix} \bigr) }{\longrightarrow}  
\overline{\cL(\partial\sigma)}\tensor\overline{\cL(\partial\delta)}\to 0.\]
From this we get 
\[ \overline{\cM(\partial(\sigma\times\delta))} \isom
\overline{\cL(\partial\sigma)}\tensor\overline{\cL(\partial\delta)}
\isom \overline{\cL}_\sigma\tensor\overline{\cL}_\delta\]
and the restriction map from $\overline{\cM}_{\sigma\times
  \delta}\isom\overline{\cL}_\sigma\tensor\overline{\cL}_\delta$ to
it is the identity.  

Now let us prove that tensoring with $\RR$ keeps the sequence above
exact. Following a similar argument in \cite{BBFK}, we show that the map
\[ \Tor^{A\otimes B}_1 (\cL(\partial\sigma)\tensor
\cL_\delta \oplus \cL_\sigma \tensor \cL(\partial\delta),\RR) \longrightarrow 
\Tor^{A\otimes B}_1
(\cL(\partial\sigma)\tensor\cL(\partial\delta),\RR) \]
is an isomorphism. For any cone $\tau$, let $\cK_\tau=
\cL([\tau],\partial\tau)$ be the kernel of the 
restriction map $\cL_\tau\to \cL(\partial\tau)$. Since $[\tau]$ is
quasi-convex in its span, $\cK_\tau$ is a free
$\cA_\tau$-module. Hence the sequence  
\[  0\longrightarrow\cK_\tau\longrightarrow \cL_\tau \]
is a free resolution of the  $\cA_\tau$-module
$\cL(\partial\tau)$. Also, the induced map $\overline{\cK}_\tau\to
\overline{\cL}_\tau$ is zero because $\overline{\cL}_\tau \isom
\overline{\cL(\partial\tau)}$. Using such  resolutions, the map of the
$\Tor^{A\otimes B}_1$ can be expressed as the isomorphism
\[ \overline{\cK}_\sigma\tensor \overline{\cL}_\delta \oplus
\overline{\cL}_\sigma \tensor \overline{\cK}_\delta
\xrightarrow{ \bigl( \begin{smallmatrix} 
    \Id & \\  & -\Id \end{smallmatrix} \bigr) } 
\overline{\cK}_\sigma\tensor \overline{\cL}_\delta \oplus
\overline{\cL}_\sigma \tensor \overline{\cK}_\delta.\] \qed

\begin{corollary} Let $\Sigma$ be a quasi-convex fan in $V$ and
$\Delta$ a quasi-convex fan in $W$. Then $\Sigma\times\Delta$ is also
  quasi-convex and there exists an isomorphism
  of $\Sym (V\times W)^* \isom \Sym 
V^*\otimes \Sym W^*$ modules
\[ \cL(\Sigma,\partial\Sigma) \tensor \cL(\Delta,\partial\Delta) \isom
\cL(\Sigma\times\Delta,\partial(\Sigma\times\Delta)).\]
\end{corollary}

{\bf Proof.}
Theorem~3.8 in \cite{BBFK} shows that a fan $\Theta$ is quasi-convex if and
only if the module $\cL(\Theta)$ is free. Since $\cL(\Sigma\times\Delta)
\isom \cL(\Sigma)\otimes\cL(\Delta)$ is free, the product fan is
quasi-convex. 

Again by \cite{BBFK} a quasi-convex $n$-dimensional fan $\Theta$
satisfies the Poincar\'e duality 
\[ \dim IH^{n-k}(\Theta) = \dim
IH^{n+k}(\Theta,\partial\Theta), \qquad k\in\ZZ.\]
Applying this to each of the three spaces in the inclusion
\[ \cL(\Sigma,\partial\Sigma) \otimes \cL(\Delta,\partial\Delta)
\hookrightarrow 
\cL(\Sigma\times\Delta,\partial(\Sigma\times\Delta)), \]
we get that both sides are free modules of the same rank (in each
degree). Thus the inclusion map must be an isomorphism. \qed

Let us assume further that $\Sigma$ and $\Delta$ are complete
projective fans with strictly convex conewise linear functions
$l_\Sigma$ and 
$l_\Delta$ on them. Then the function $l=\pi_1^*l_\Sigma + \pi_2^*
l_\Delta $ is strictly convex on the product fan
$\Sigma\times\Delta$. We wish to state that if $Q_{l_\Sigma}$ and
$Q_{l_\Delta}$ both satisfy the \HRM bilinear relations
on the respective cohomology spaces, then so does $Q_l$. For this
we need to assume that the representation of sections
$\cL(\Sigma\times\Delta)\isom \cL(\Sigma)\otimes\cL(\Delta)$ as
conewise polynomial functions on some subdivision is induced by the
representation of sections $\cL(\Sigma)$ and $\cL(\Delta)$ as
functions, and that the evaluation map satisfies:
\[ \langle h_1\otimes h_2\rangle_{\Sigma\times\Delta} = \langle h_1\rangle_\Sigma\cdot \langle h_2\rangle_\Delta,
\qquad h_1\in IH(\Sigma), \; h_2\in IH(\Delta).\]
If these conditions are satisfied then we have:

\begin{corollary}\label{cor-prod} If $Q_{l_\Sigma}$ and
$Q_{l_\Delta}$ both satisfy the \HRM bilinear relations
then so does $Q_l$.
\end{corollary}

{\bf Proof.} We choose a homogeneous basis $\{p_i\}$ of $IP(\Sigma)$,
orthogonal 
with respect to the bilinear form $B_{l_\Sigma}(x,y) = \langle l_\Sigma^k x\cdot
y\rangle$, and a similar orthogonal basis $\{q_j\}$ of $IP(\Delta)$. Then it
suffice to show that $Q_l$ satisfies the \HRM relations on the subspace
$S\subset IH(\Sigma\times\Delta)$ with basis 
\[ \{ l_\Sigma^\alpha p\otimes l_\Delta^\beta q \}_{0\leq \alpha\leq
  k_1, \; 0\leq \beta\leq k_2} \]
generated by one $p= p_{i_0}$ and one $q=q_{j_0}$. Since $l$ acts by
multiplication with $l_\Sigma\otimes 1 + 1\otimes l_\Delta$, it is
clear that the signature of $Q_l$ on $S$ only depends on the numbers
$k_1, k_2$ and the signs of $Q_{l_\Sigma}(p), Q_{l_\Delta}(q)$. 

One can give a geometric proof that $Q_l$ satisfies the \HRM relations on
$S$ by comparing it with the quadratic form $Q_A$ defined on the
singular cohomology $H(\CC\PP^{k_1}\times \CC\PP^{k_2},\CC)$, where the
Lefschetz action is given by intersection with an ample divisor $A$, and
the pairing is the usual cup product. Then $Q_A$ satisfies the
classical Hodge-Riemann bilinear relations \cite{GH}. The
combinatorial analog of this is to construct two simplicial fans
$\Pi^{k_1}$ and $\Pi^{k_2}$
corresponding to the toric varieties $\CC\PP^{k_1}$ and
$\CC\PP^{k_2}$, together with strictly convex functions on them. Then
the \HRM relations \cite{M} for the simplicial fan
$\Pi^{k_1} \times \Pi^{k_2}$ imply the the same relations for $Q_l$
above. \qed

\subsection{Skew products}

Consider two fans $\Sigma$ in $V$ and $\Delta$ in $W$. Let
$\phi:|\Sigma|\to W$ be a map which is linear on each cone
$\sigma\in\Sigma$, and let $\Gamma_\sigma(\phi)$ be the graph of
$\phi$ restricted to $\sigma$. We define the fan $\Sigma\times_\phi
\Delta$ in $V\times W$: 
\[ \Sigma\times_\phi \Delta = \{ \Gamma_\sigma(\phi)+\delta\mid
\sigma\in\Sigma, \; \delta\in\Delta\}.\]
Then if $\phi$ is the zero map, we recover the usual product. For a 
general $\phi$ we only get a conewise linear isomorphism 
\[ \Phi: \Sigma\times\Delta \longrightarrow \Sigma\times_\phi \Delta \]
defined by $|\Phi|(v,w) = (v,w+\phi(v))$ for $v\in |\Sigma|$ and
$w\in|\Delta|$. 

Since $\Phi$ is a conewise linear isomorphism of fans, it induces an
isomorphism of vector spaces
\[ \cL( \Sigma\times_\phi \Delta) \isom
\cL(\Sigma\times\Delta) \isom\cL(\Sigma)\otimes\cL(\Delta).\] 
Let $A$ and $B$ be the rings of polynomial functions on $V$ and $W$,
respectively. Then the isomorphism of global sections is in fact an
isomorphism of $A$-modules. A functions $f\in B$ acts on the tensor
product by multiplication with $\Phi^*(1\otimes f)$. In particular, if
$f\in B$ is linear, then  
\[ \Phi^*(1\otimes f) = 1\otimes f + f\circ\phi\otimes1. \]

\subsection{Local product structures}

We say that a fan $\Sigma$ in $V$ has a local product structure at a cone
$\sigma\in\Sigma$ if 
\[ \overline{\Star}_\Sigma(\sigma)\isom \Lambda\times_\phi [\sigma]\]
for some fan $\Lambda$ in $W$, where $V=W\oplus\Span(\sigma)$,
and some function $\phi:|\Lambda|\to \Span(\sigma)$. Equivalently,
$\Sigma$ has a local product 
structure at $\sigma$ if every cone $\tau\in\Sigma$ containing
$\sigma$ can be written as a direct sum $\tau=\rho\oplus\sigma$ for some
$\rho\in\Link_\Sigma(\sigma)$. If $\pi:
\overline{\Star}_\Sigma(\sigma)\to \Lambda$ is the projection from
$\Span(\sigma)$ then $\pi$ restricts to a conewise linear
isomorphism $\Link_\Sigma(\sigma) \to \Lambda$.

As an example, when we perform a star subdivision of the fan
$\Sigma$ at $\RRp v$, where $v$ lies in the relative interior of a
cone $\sigma$, we require $\Sigma$ to have a local product structure
at $\sigma$. The subdivision $\hat{\Sigma}$ then also has a local
product structure at $\RRp v\in \hat{\Sigma}$.

\begin{lemma} \label{lem-loc-prod} Assume that $\Sigma$ has a local
  product structure at a simplicial cone $\sigma\in\Sigma$. Let $\Delta =
  \overline{\Star}_\Sigma(\sigma)$; then with notation as above, 
\[ \pi^* \cL_\Lambda \isom \cL_\Delta. \]
Moreover, $\pi^*$ induces an isomorphism of $A$-modules
\[ \cL(\Delta) \isom \cL(\Lambda)\otimes_B A,\]
hence also an isomorphism $IH(\Delta) \isom IH(\Lambda)$.
\end{lemma}

{\bf Proof.} A cone $\tau\in\Delta$ has the form $\tau=
\Gamma_\rho(\phi)\oplus \sigma'$ for some $\rho\in\Lambda$ and $\sigma'$ a
face of $\sigma$. Since $\sigma'$ is simplicial,
$\cL_{\sigma'}\isom\cA_{\sigma'}$ and we have an isomorphism of
$\cA_\tau$ modules
\[\cL_\tau \isom \cL_{\Gamma_\rho(\phi)}\otimes \cA_{\sigma'} \isom
\cL_\rho\otimes_{\cA_\rho} \cA_{\tau}.\]
This proves the first isomorphism. The isomorphism above also induces
an isomorphism of $A$-modules
\[ \cL_\tau \isom \cL_\rho\otimes_B A.\]

Let $M\subset\Lambda$ be the set of maximal cones, i.e., cones which
are not proper faces of any other cones. Then the cones
$\Gamma_\rho(\phi)\oplus 
\sigma$ for $\rho\in M$ are the maximal cones of $\Delta$. Using the
covering by maximal cones, $\cL(\Delta)$ is the kernel of the
restriction map
\[ \bigoplus_{\rho\in M} \cL_{\Gamma_\rho(\phi)\oplus \sigma} \longrightarrow
\bigoplus_{\rho_1,\rho_2\in M} \cL_{\Gamma_{\rho_1}(\phi)\oplus \sigma \cap
  \Gamma_{\rho_2}(\phi)\oplus \sigma}. \]
We can write this sequence of $A$-modules as
\[ (\bigoplus_{\rho\in M} \cL_\rho \longrightarrow
\bigoplus_{\rho_1,\rho_2\in M} \cL_{\rho_1\cap\rho_2})\otimes_B A. \]
Since the kernel of the map in the parentheses is the $B$-module
$\cL(\Lambda)$, we get the stated isomorphism of global sections and
intersection cohomologies. \qed

\section{Distinguished embeddings}

The Poincar\'e pairing, and hence the definition of
the quadratic form $Q_l$, depend on the representation of sections of
$\cL_\Sigma$ as conewise polynomial functions on some subdivision of
$\Sigma$. The goal of this section is to fix one such representation.

\subsection{Primitive embeddings} By the definition of the sheaf
$\cL$, for any cone $\sigma$ we have 
\[ \cL_\sigma \isom \overline{\cL(\partial\sigma)}\otimes\cA_\sigma.\]
We start by recalling how $\overline{\cL(\partial\sigma)}$ can be
identified with the primitive cohomology of the ``flattened boundary
fan'' \cite{BBFK, BL}.

By restricting to the span of $\sigma$ if necessary, we may assume
that $\sigma$ is an $n$-dimensional cone in $V$. Let $v$ 
be a vector in the relative interior of $\sigma$, $W= V/\RR v$, and
let $\pi:V\to W$ be the projection. Then $\pi$ induces a conewise linear
isomorphism    
\[ \pi: \partial\sigma \to \Lambda\]
for a complete fan $\Lambda$ in W. If we let $B = \Sym W^*$ be the
ring of polynomial functions on $W$ and $A=B[x]$ for some linear
function $x$ not vanishing on $v$, then $\pi^*$ defines an isomorphism
of $B$-modules
\[ \cL(\partial\sigma) \isom \cL(\Lambda).\]
Multiplication of sections $\cL(\partial\sigma)$ with $x$
corresponds via this isomorphism to multiplication of sections
$\cL(\Lambda)$ with a strictly convex conewise linear function
$\lambda$ on $\Lambda$, where
\[ \lambda = x \circ (\pi|_{\partial\sigma})^{-1}.\]
Thus, if we assume the hard Lefschetz theorem for the lower dimensional fan
$\Lambda$ (an assumption we will make everywhere in this section), we have
\[ \overline{\cL(\partial\sigma)} \isom
\overline{\cL(\Lambda)}/\lambda\cdot \overline{\cL(\Lambda)} \isom
IP_\lambda(\Lambda).\] 
We note that this construction does not depend on a particular choice
of $x$; indeed, modulo the maximal ideal $I_B\subset B$ generated by
linear functions, $x$ is determined up to a constant factor.

Now let $[\hat\sigma]$ be the star subdivision of $[\sigma]$ at $\RRp
v$. Then by Lemma~\ref{lem-loc-prod}, the pullback $\pi^*$ induces an
isomorphism  
\[ \overline{\cL([\hat\sigma])} \isom \overline{\cL(\Lambda)} =
IH(\Lambda),\]
hence multiplication with $\pi^*(\lambda)$ defines a \LO on
$\overline{\cL([\hat\sigma])}$ centered at degree $\dim\sigma - 1$. 

\begin{definition}\label{def-prim-emb} Let $s:\hat{\Sigma}\to\Sigma$
  be the star 
  subdivision of a fan $\Sigma$ at $\RRp v$, where $v$ lies in the
  relative interior of a cone $\sigma\in\Sigma$. We say that a direct
  embedding (i.e, an embedding as a direct summand of locally free
  $\cA_\Sigma$ modules)  
\[\cL_\Sigma \subset s_* \cL_{\hat{\Sigma}}\]
is a {\em primitive embedding} if the induced embedding
\[ \overline{\cL([\sigma])} \subset \overline{\cL([\hat\sigma])} \]
identifies $\overline{\cL([\sigma])}$ with the $\lambda$-primitive part
of $\overline{\cL([\hat\sigma])}$. Here $[\hat\sigma]\subset\hat\Sigma$
is the star subdivision of $[\sigma]$ at $\RRp v$, and $\lambda$ is a
conewise linear function on $[\hat\sigma]$ such that $\lambda(v)=0$, 
$\lambda|_{\partial\sigma} = x|_{\partial\sigma}$ for some linear
function $x\in A$, $x(v)\neq 0$. 
\end{definition}

\subsection{Distinguished subdivisions} \label{sec-dist}
We now specify certain distinguished simplicial subdivisions of a
fan $\Sigma$ and distinguished representations of sections of
$\cL_\Sigma$ as conewise polynomial functions on the subdivision.

\begin{definition} A {\em distinguished subdivision} of a fan $\Sigma$ 
is a sequence
\[  \Sigma_0 \stackrel{s_1}{\longrightarrow} \Sigma_1 \stackrel{s_2}{\longrightarrow} \ldots
\stackrel{s_N}{\longrightarrow} \Sigma_N = \Sigma,\]
where $s_i:\Sigma_{i-1}\to\Sigma_{i}$ is a star subdivision at $\RRp
v_i$ for $i=1,\ldots,N$, satisfying the following conditions
\begin{enumerate}
\item $\Sigma_0$ is simplicial.
\item There exists a subfan $\Sigma^s\subset\Sigma$, the ``singular''
  subfan, such that the sequence of subdivisions $s_i$ is induced by a
  barycentric subdivision of $\Sigma^s$. In other words, $v_i\in
  |\Sigma^s|$, and the sequence of star subdivisions of $\Sigma^s$ at
  $\RRp v_i$ for $i=N,\ldots,1$ is a 
  barycentric subdivision of $\Sigma^s$.
\item The subfan $\Sigma^s$ has the property that $|\Sigma^s|\cap
  \sigma$ is a face of $\sigma$ for any cone $\sigma\in\Sigma$.
\end{enumerate}
\end{definition}

The second condition of the definition implies that any cone $\tau\in
\Sigma^s$ of 
dimension at least $2$ contains a unique $v_i$ in its relative
interior for some $1\leq i\leq N$. We denote this $v_i$ by $v_\tau$. 

As an example, a barycentric subdivision  is a distinguished
subdivision of $\Sigma$, with $\Sigma^s = \Sigma$. Also, if 
$s = s_{N}\circ\ldots \circ s_1$ is a distinguished subdivision of
$\Sigma=\Sigma_N$, then $t = s_{N-1}\circ\ldots \circ s_1$ is a
distinguished subdivision of $\Sigma_{N-1}$, with $\Sigma_{N-1}^s =
\Sigma^s- \{\sigma_{N}\}$, where $\sigma_{N}\in \Sigma$ is
the cone containing $v_{N}$ in its relative interior. 

Suppose we have a distinguished subdivision of $\Sigma$. Given a cone
$\sigma\in\Sigma$, we can express $\sigma = 
\tau\oplus\rho$, where $\tau=\sigma\cap|\Sigma^s|$ and $\rho$ is a
simplicial cone. Using the K\"unneth theorem and the definition of the
sheaf $\cL$, we have
\[ \cL_\sigma \isom \cL_\tau\otimes\cL_\rho \isom
\overline{\cL(\partial\tau)} \otimes_{\cA_\tau} \cA_\sigma.\]
Thus, a set of $\cA_\tau$ module generators of $\cL(\partial\tau)$ can
be chosen as the $\cA_\sigma$ module generators of $\cL_\sigma$. When
representing sections of $\cL$ as functions, we impose this condition
on generators explicitly.

\begin{definition} Let $s = s_{N}\circ\ldots \circ s_1$ be a
distinguished subdivision of a fan $\Sigma$. Given an embedding of
$\cA_\Sigma$ modules 
\[ \cL_\Sigma \subset s_* \cA_{\Sigma_0},\]
we have a representation of sections of $\cL_\Sigma$ as conewise
polynomial functions on $\Sigma_0$. We say that this representation is
{\em distinguished} if it satisfies the following condition. For any
cone $\sigma\in\Sigma$, write $\sigma= \tau\oplus\rho$, where $\tau =
\sigma\cap|\Sigma^s|$ and $\rho$ is a simplicial cone. Then the
condition is that there exist generators of the $\cA_\sigma$ module
$\cL_\sigma$ represented by functions on $\partial\tau$, pulled back
to $\sigma$ by the projection map 
\[ \pi: \Span \sigma \longrightarrow \Span\sigma/ (\RR v_\tau + \Span\rho).\]
\end{definition}

We remark that choosing a distinguished embedding $\cL_\Sigma \subset
s_* \cA_{\Sigma_0}$ is equivalent to choosing such an embedding for
the sheaves restricted to the subfan $\Sigma^s$.

\begin{lemma}\label{lem-loc-prod1} Consider a distinguished subdivision $s =
s_{N}\circ\ldots \circ s_1$ of a fan $\Sigma$ and a distinguished
representation of sections of $\cL_\Sigma$ as functions on
$\Sigma_0$. Let $\rho\in\Sigma$ be a cone such that $\rho\cap\Sigma^s
= \{0\}$, and denote $\Delta = \overline{\Star}_\Sigma \rho$. Then
$\Sigma$ has a local product structure at $\rho$:
\[ \Delta = \Lambda\times_\phi [\rho], \]
and there exist generators of the $A$ module $\cL(\Delta)$ represented
by functions on (the subdivision of) $\Lambda$, pulled back to
$\Delta$ by the projection 
\[ \pi: \Delta\longrightarrow\Lambda .\]
\end{lemma}

{\bf Proof.} Since any cone $\sigma\in\Sigma$ containing $\rho$ can be
written as $\sigma=\rho\oplus\tau$, it follows that $\Sigma$ has a
local product structure at $\rho$. Moreover, since the representation
of sections of $\cL_\Sigma$ is distinguished, generators of
$\cL_\sigma$ can be chosen as functions on $\tau$ pulled back by the
projection $\sigma\to\tau$ from $\Span(\rho)$. Thus,
\[ \cL_\sigma = \cL_\tau \otimes_B A,\]
where $B$ is the space of polynomial functions on $V/\Span(\rho)$ and
$\cL_\tau$ is a $B$-module of functions on $\tau$. Letting $M$ be the
set of maximal cones in $\Link_\Sigma\rho$, we get as in the proof of
Lemma~\ref{lem-loc-prod} that $\cL(\Delta)$ is the kernel of the
map 
\[ (\bigoplus_{\tau\in M} \cL_\tau \longrightarrow\bigoplus_{\tau_1,\tau_2\in M}
\cL_{\tau_1\cap\tau_2})\otimes_B A, \] 
hence a set of $B$-module generators of the kernel of the map in the
parentheses forms a set of  generators for the $A$-module $\cL(\Delta)$. \qed

\begin{lemma}\label{lem-prim-emb} Let $s = s_{N}\circ\ldots \circ s_1$
  be a distinguished 
subdivision of a fan $\Sigma$. Assume that we have an embedding 
\[ \cL_{\Sigma_{N-1}} \subset (s_{N-1}\circ\ldots \circ s_1)_*
\cA_{\Sigma_0} \] 
giving a distinguished representation of sections of $\cL_{\Sigma_{N-1}}$ as
conewise polynomial functions on the fan $\Sigma_0$. Then there exists
a primitive embedding
\[ \cL_\Sigma \subset (s_N)_* \cL_{\Sigma_{N-1}},\]
such that the induced representation of sections of $\cL_\Sigma$ as
functions on $\Sigma_0$ is distinguished.
\end{lemma}

{\bf Proof.} Note that to construct an embedding $\cL_\Sigma \subset
(s_N)_* \cL_{\Sigma_{N-1}}$, we have to choose a subsheaf of $(s_N)_*
\cL_{\Sigma_{N-1}}$ satisfying the conditions for being an equivariant
intersection cohomology sheaf. This subsheaf is then isomorphic to
$\cL_\Sigma$. Also, it suffices to construct the embedding of the
sheaves when restricted to the subfan $\Sigma^s$, so we may assume that
$\Sigma=\Sigma^s$ and $s_N: \Sigma_{N-1}\to\Sigma$ is the star
subdivision at $\RRp v$, where $v$ lies in the relative interior of a
maximal cone $\sigma$. For any cone $\tau\in\Sigma$, $\tau\neq\sigma$,
we set $\cL_\tau = ((s_N)_* \cL_{\Sigma_{N-1}})_\tau$. If
$\tau=\sigma$, then by the previous lemma, generators of $((s_N)_*
\cL_{\Sigma_{N-1}})_\sigma = \cL([\hat\sigma])$ can be chosen as
functions on $\partial\sigma$, pulled back by the projection map
$\pi:\Span\sigma\to \Span\sigma/\RR v$. We choose a subset of these
generators, representing a basis for the primitive part of
$\overline{\cL([\hat\sigma])}$ with respect to a conewise linear
function $\lambda$ as in Definition~\ref{def-prim-emb}. Such a subset
generates the module $\cL_\sigma$. \qed

\begin{definition} Consider the following data associated to a fan
  $\Sigma$:
\begin{enumerate}
\item[(a)] A distinguished subdivision 
\[  \Sigma_0 \stackrel{s_1}{\longrightarrow} \Sigma_1 \stackrel{s_2}{\longrightarrow} \ldots
\stackrel{s_N}{\longrightarrow} \Sigma_N = \Sigma\]
of the fan $\Sigma$.
\item[(b)] A sequence of primitive embeddings 
\[ e_1: \cL_{\Sigma_1} \subset (s_1)_* \cA_{\Sigma_0}, \qquad e_i:
\cL_{\Sigma_i} \subset (s_i)_* \cL_{\Sigma_{i-1}}, \quad i=2,\ldots,N,\] 
such that for any $i=1,\ldots,N$, the induced representation of
sections of $\cL_{\Sigma_i}$ as conewise polynomial functions on the
fan $\Sigma_0$ is distinguished.
\end{enumerate}
We call such data $(\{s_i\},\{e_i\})$ a {\em distinguished pair} on
$\Sigma$.
\end{definition}

The following properties of distinguished pairs are straight-forward:
\begin{enumerate}
\item Every fan $\Sigma$ has a distinguished pair. It suffices to take
a barycentric subdivision of $\Sigma$ and construct the primitive
embeddings by induction on $i$ as in Lemma~\ref{lem-prim-emb}. 
\item Given a distinguished pair on $\Sigma=\Sigma_N$, we get a distinguished
  pair on $\Sigma_i$ for any $i=1,\ldots,N-1$.
\item A distinguished pair on $\Sigma$ can be restricted to a
  distinguished pair on a subfan $\Delta\subset\Sigma$. We take
  $\Delta^s = \Delta\cap \Sigma^s$ and restrict the subdivisions
  and the sheaves.
\item If $\Phi: \Delta\to\Sigma$ is a conewise linear isomorphism,
  then a distinguished pair on $\Sigma$ can be pulled back to a
  distinguished pair on $\Delta$. 
\end{enumerate}

As an application of these properties, consider a fan $\Sigma$ with a
distinguished pair on it and a cone $\rho\in\Sigma$ such that
$\rho\cap\Sigma^s = \{0\}$. Then $\Sigma$ has a local product
structure at $\rho$:
\[ \overline{\Star}_\Sigma \rho  = \Lambda\times_\phi [\rho]. \]
We claim that the distinguished pair on $\Sigma$ induces a
distinguished pair on $\Lambda$. Indeed, by restriction we get a
distinguished pair on $\Link_\Sigma \rho \subset\Sigma$, and then we
pull this back to $\Lambda$ by the conewise linear isomorphism 
$\Lambda\to\Link_\Sigma \rho$.

\subsection{Poincar\'e duality}

Let $\Sigma$ be a quasi-convex fan. The Poincar\'e duality pairing 
\[ IH^{n-k}(\Sigma)\times IH^{n+k}(\Sigma,\partial\Sigma)\longrightarrow\RR\]
in \cite{BBFK} is constructed by representing sections of $\cL(\Sigma)$ as
conewise polynomial functions on a barycentric subdivision of
$\Sigma$, and then multiplying these functions. Now suppose we have a
barycentric subdivision of $\Sigma$:
\[  \Sigma_0 \stackrel{s_0}{\longrightarrow} \Sigma_1 \stackrel{s_1}{\longrightarrow} \ldots
\longrightarrow \ \Sigma_{N-1} \stackrel{s_{N-1}}{\longrightarrow} \Sigma_N = \Sigma.\]
If we follow \cite{BBFK}, to construct the Poincar\'e duality paring
on the fan $\Sigma_{N-1}$, we would have to consider its barycentric
subdivision. Instead, we want to represent sections of
$\cL_{\Sigma_{N-1}}$ as conewise polynomial functions on the fan
$\Sigma_0$. The proof of Poincar\'e duality in \cite{BBFK}
also works in this case. We indicate only what modifications have to
be made.

Let us assume that we have a distinguished pair on $\Sigma$.
Thus we represent sections of $\cL_\Sigma$ as conewise polynomial
functions on 
$\Sigma_0$. We claim that multiplication of these functions followed
by the evaluation map $IH^n(\Sigma_0)\to\RR$ induces a nondegenerate
pairing 
\[ IH^{n-k}(\Sigma)\times IH^{n+k}(\Sigma,\partial\Sigma)\longrightarrow\RR\]
for any $k\in\ZZ$.

The proof in \cite{BBFK} has two parts. First, assuming that the
statement holds for fans $[\sigma]$ where $\sigma$ is an
$n$-dimensional cone, it is shown to hold for any quasi-convex fan of
dimension 
$n$. This part of the proof also works in our case because it does not
depend on the representation of sections of $\cL_\Sigma$ as
functions. The second part of the proof deduces Poincar\'e
duality for fans $[\sigma]$, assuming that it holds for fans of
dimension less than $n$. Here is where the representation of sections
of $\cL$ as functions is important and the proof needs to be modified.

Assume that $\Sigma = [\sigma]$ for a cone $\sigma$. Then $\sigma =
\tau\oplus\rho$, where $\tau=\sigma\cap|\Sigma^s|$ and $\rho$ is a
simplicial cone. If $\tau=\sigma$, then the sequence of subdivisions
is a barycentric 
subdivision of $[\sigma]$ and the proof of \cite{BBFK} applies. If
$\tau\neq\sigma$ then since the representation of sections is
distinguished, we have the identifications
\begin{align*} \cL_\sigma &= \cL_\tau\otimes\cA_\rho,\\
\cL([\sigma],\partial\sigma) &=
\cL([\tau],\partial\tau)\otimes\cA([\rho],\partial\rho).
\end{align*}
By induction on the dimension we may
assume that the Poincar\'e duality pairing is nondegenerate on $[\tau]$
and $[\rho]$, hence it is nondegenerate on $[\sigma]=[\tau]\times[\rho]$.

\section{Restriction to facets}

In \cite{M} McMullen proved that for simplicial projective fans the
\HRM relations in dimension $n-1$ imply the Hard Lefschetz theorem in
dimension $n$. The proof, written in the language of polytope
algebras, relies on restricting an element of the polytope algebra to
the facets of a polytope. We explain this idea in terms of fans and
then give a partial generalization to nonsimplicial fans.

We start with a more general situation. Let $\Sigma$ be a complete  
fan with a distinguished pair. Let $\RRp v\in\Sigma$ be a
1-dimensional cone such that $\RRp v\notin\Sigma^s$, and let $\Delta =
\overline{\Star}_\Sigma(\RRp v)$. Then $\Sigma$ has a local product
structure at $\RRp v$:
\[ \Delta = \Lambda\times_\phi[\RRp v].\]
We let $\pi:\Delta\to\Lambda$ be the projection, inducing an
isomorphism
\[ IH(\Delta) \isom IH(\Lambda).\]
For $h\in IH(\Sigma)$, we denote by $h|_\Lambda$ the image of $h$
under the restriction map $IH(\Sigma)\to IH(\Delta)\isom
IH(\Lambda)$. In terms of sections, the restriction of $h$ to $\Delta$
can be represented by a function of the form $\pi^*(f)$ (see
Lemma~\ref{lem-loc-prod1}); then $f$ represents $h|_\Lambda$. Note
also that the distinguished pair on $\Sigma$ induces a distinguished
pair on $\Lambda$, thus we get the non-degenerate Poincar\'e pairing,
and by induction on dimension, the \HRM relations associated with a
strictly convex conewise linear function  on $\Lambda$. 

Let $\lambda \in \cA^2(\Sigma)$ be a conewise linear function supported
on $\Delta$, such that $\lambda(v)>0$. Multiplication with $\lambda$
defines a degree $2$ map $IH(\Sigma)\to IH(\Sigma)$, with image lying
in $IH(\Delta,\partial\Delta)$. 

\begin{lemma}\label{lem-reduct} For any $h,g\in IH(\Sigma)$ we have 
\[ \langle \lambda h \cdot g\rangle_\Sigma = \langle h|_\Lambda\cdot g|_\Lambda\rangle_\Lambda\]
up to a positive constant factor.
\end{lemma}

{\bf Proof.} By assumption, we have a simplicial subdivision $\Sigma_0$
of $\Sigma$ such that cohomology classes are represented by conewise
polynomial functions on $\Sigma_0$. Recall the construction of the
evaluation map in the top degree cohomology. When constructing the functions
$\Phi_\sigma$ for $\sigma\in\Sigma_0$ an $n$-dimensional cone such
that $\RRp v\subset \sigma$, we set $H_{\sigma,1} =
\lambda|_\sigma$. Then for any function $f$ on $\Lambda$ we have
\[ \langle \lambda \pi^*(f)\rangle_\Sigma = \langle f\rangle_\Lambda\]
up to a constant factor which depends on the metrics on $\Lambda^n V^*$
and $\Lambda^{n-1} (V/\RR v)^*$. \qed

We will assume that the constant in the previous lemma is always $1$.
For later use we also prove 

\begin{lemma}\label{lem-deg2-isom} Multiplication with $\lambda$
  defines a degree $2$ 
  linear map $\lambda: IH(\Sigma)\to IH(\Delta,\partial\Delta)$. This
  map factors through the restriction map $IH(\Sigma)\to IH(\Delta)$
  and the map $\lambda|_\Delta: IH(\Delta)\to
  IH(\Delta,\partial\Delta)$ is a degree $2$ isomorphism.
\end{lemma}

{\bf Proof.} Since $\lambda$ has support on $\Delta$, the only
statements that needs a proof is that $\lambda|_\Delta$ defines an
isomorphism. Consider the bilinear form on $IH(\Delta)$
\[ (h,g) \longrightarrow \langle \lambda|_\Delta h\cdot g\rangle_\Sigma = \langle h|_\Lambda\cdot
g|_\Lambda\rangle_\Lambda.\] 
The right hand side is the Poincar\'e pairing on $\Lambda$, hence the
map  $\lambda|_\Delta: IH(\Delta)\to  IH(\Delta,\partial\Delta)$ is
injective and the image can be identified with the dual of
$IH(\Delta)$. By Poincar\'e duality applied to the fan $\Delta$,
$IH(\Delta,\partial\Delta)$ is the dual of $IH(\Delta)$, hence the map
is an isomorphism. \qed

Let us now return to McMullen's argument. Assume that $\Sigma$ is a
simplicial fan, and let $l\in\cA^2(\Sigma)$ be a strictly convex
function, positive on 
$V-\{0\}$. For each $1$-dimensional cone $\RRp v_i\in\Sigma$ we let
$\lambda_i\in\cA^2(\Sigma)$ be supported on
$\Delta_i=\overline{\Star}_\Sigma(\RRp v_i)$ 
so that 
\[ l = \sum_i \lambda_i.\]
Associated with $\Delta_i$ we have the fan $\Lambda_i$ and the
projection map $\pi_i$ as before.

Suppose there exists an $h\in IH^{n-k}(\Sigma)$ such that
\[ l^k\cdot h = 0.\] 
We restrict this equality to $\Lambda_i$. Note that 
we may change $l$ by a global linear function so that $l|_{\Delta_i} =
\pi_i^*(l_i)$ for some strictly convex function $l_i$ on
$\Lambda_i$. It follows that 
\[ h|_{\Lambda_i} \in IP_{l_i}(\Lambda_i).\]
If we assume the \HRM relations in dimension $n-1$, we get 
\[ (-1)^\frac{n-k}{2} \langle l_i^{k-1} h|_{\Lambda_i} \cdot
h|_{\Lambda_i}\rangle_{\Lambda_i} =  (-1)^\frac{n-k}{2}
\langle \lambda_i l^{k-1} h \cdot h\rangle_\Sigma \geq 0,\]
and the equality holds if and only if $h|_{\Lambda_i}=0$. Summing over all $i$,
\[ \sum_i  (-1)^\frac{n-k}{2} \langle \lambda_i l^{k-1} h \cdot h\rangle_\Sigma =
(-1)^\frac{n-k}{2} \langle  l^{k} h \cdot h\rangle_\Sigma.\] 
Since $l^k h = 0$, we get that $h|_{\Lambda_i} = 0$ for all $i$. 

To finish McMullen's proof, it suffices to show that if $h|_{\Lambda_i}
= 0$ for all $i$ then $h=0$. This can be done by showing that
$\cA(\Sigma)$ is generated as an $A$ module by elements having support
on $\Delta_i$; then the Poincar\'e pairing between $h$ and any such
element is zero. An explicit set of generators supported on $\Delta_i$
was constructed by Timorin in \cite{T1}.

In the nonsimplicial case we have the following partial
generalization of this argument. If $\Sigma$ has a distinguished pair
on it and $\RRp v_i\in\Sigma$ is a
$1$-dimensional cone such that  $\RRp v_i\notin\Sigma^s$, then
$\Sigma$ has a local product structure at $\RRp v_i$. We use the same
notation $\Delta_i, \Lambda_i, \pi_i, l_i$ associated with the ray $\RRp
v_i$ as in the simplicial case. 

\begin{lemma}\label{lem-rest} Let $\Sigma$ be a complete
  $n$-dimensional fan with a 
  distinguished pair 
  on it, such that $\Sigma^s\subset [\sigma]$ for a single cone
  $\sigma\in\Sigma$, and let $l\in\cA^2(\Sigma)$ be a strictly convex
  function. Then the \HRM relations in dimension $n-1$,
  imply that any $h\in IH^{n-k}(\Sigma)$ such that $l^k h=0$ must
  satisfy
\[ h|_{\Lambda_i} = 0 \]
for all $i$ such that $\RRp v_i \in \Sigma - [\sigma]$.
\end{lemma}

{\bf Proof.} Let $\RRp v_i \in \Sigma - [\sigma]$ be a $1$-dimensional
cone. Restricting the equality $l^k h=0$ to $\Lambda_i$, we get that
$h|_{\Lambda_i} \in IP_{l_i}(\Lambda_i)$. 

Now we modify the function $l$ by a global linear function 
so that $l|_\sigma = 0$, $l(x)>0$ for $x\notin\sigma$, and write
\[ l = \sum_{\RRp v_i\in\Sigma- [\sigma]}\lambda_i \]
with $\lambda_i$ as before. Then
\begin{gather*} 0 = (-1)^\frac{n-k}{2} \langle  l^{k} h \cdot h\rangle_\Sigma = 
\sum_{\RRp v_i\in\Sigma- [\sigma]}  (-1)^\frac{n-k}{2} \langle \lambda_i
l^{k-1} h \cdot h\rangle_\Sigma \\
= \sum_{\RRp v_i\in\Sigma- [\sigma]}
\langle l_i^{k-1} h|_{\Lambda_i} \cdot h|_{\Lambda_i}\rangle_{\Lambda_i}.
\end{gather*}
The primitivity of $h|_{\Lambda_i}$ implies that all
these restrictions must be zero. \qed

\section{Proof of the main theorem}

We now give a proof of Theorem~\ref{thm-HRM}.

Let $\Sigma$ be a complete fan with a distinguished pair, and let $l$
be a conewise linear strictly convex function on $\Sigma$.  Let 
\[  \Sigma_0 \stackrel{s_1}{\longrightarrow} \Sigma_1 \stackrel{s_2}{\longrightarrow} \ldots
\stackrel{s_N}{\longrightarrow} \Sigma_N = \Sigma\]
be the distinguished subdivision of $\Sigma$, where $s_i$ is the star
subdivision at the ray $\RRp v_i$. We choose for
$i=N-1,\ldots,0$ a strictly convex conewise linear function $l_i$ on
$\Sigma_i$ such that $l_i-l_{i+1}$  
is supported on $\Star_{\Sigma_{i}} v_{i+1}$, $l_N=l$.
Since for all $i=0,\ldots,N$  we have a distinguished pair on the fan
$\Sigma_i$, we can define the quadratic form $Q_{l_i}$ on
$IH(\Sigma_i)$. We prove by induction on $i$ that $Q_{l_i}$ satisfies
the \HRM relations on $IH(\Sigma_i)$. The case $i=0$ follows by our
assumption that the relations hold for simplicial fans. The induction
step then is to prove the  \HRM relations for $Q_{l_N} = Q_l$,
assuming that $Q_{l_{N-1}}$ satisfies these relations. We may also
assume by induction on the dimension that the   \HRM relations are
satisfied in the case of lower dimensional fans.

To simplify notation, let $\hat{\Sigma} \to \Sigma$ be a star
subdivision at $\RRp v$, where $v$ lies in the relative interior of a
cone $\sigma\in\Sigma$ ($\hat{\Sigma} = \Sigma_{N-1}$ in the previous
notation), and let $\hat{l}$ be a strictly convex conewise linear
function on $\hat{\Sigma}$ such that $l-\hat{l}$ is supported on
$\Star_{\Sigma} \sigma$. 

Let us denote  $\Delta = \overline{\Star}_\Sigma(\sigma)$ and let $\hat{\Delta}
\subset\hat{\Sigma}$ be the star subdivision of $\Delta$ at $\RRp v$.
We write 
\[ \hat\Delta = \Theta\times_\phi[\RRp v],\]
and let $\pi: \hat\Delta\to\Theta$ be the projection from $\RR v$.
Changing $l$ and  $\hat{l}$ by a linear function if necessary, we may
assume that $\hat{l}$ 
restricted to $\hat{\Delta}$ is the pullback of a strictly convex
function $l_\Theta$ on $\Theta$.
% hence it defines a \LO on
% $IH(\hat{\Delta})\isom IH(\Theta)$ centered at
% degree $n-1$. 

Next we define an auxiliary fan $\overline{\Delta}$:
\[ \overline{\Delta} = \Delta\cup \Delta',\]
where
\[ \Delta' = \{ \tau+\RRp(-v) \mid \tau\in\partial\Delta \} \cup
\partial\Delta.\]
On $\overline{\Delta}$ we have a conewise linear
function $\overline{l}$:
\[ \overline{l}(x) = \begin{cases} l(x) & \text{if $x\in |\Delta|$,} \\
l_\Theta(\pi(x)) & \text{if $x\in|\Delta'|$.} \end{cases}\]
Then $\overline{l}$ is strictly convex on $\overline{\Delta}$. Indeed,
it suffices to prove that $\overline{l}(x) > l_\Theta(\pi(x))$ for
$x$ in the interior of $|\Delta|$; this follows from the inequality
$l(x)>\hat{l}(x)$ for $x$ in the interior of $|\Delta|$. 

We claim that the distinguished pair on $\Sigma$ naturally induces
a distinguished pair on $\overline{\Delta}$. The restriction of the
distinguished pair on $\Sigma$ to $\Delta$ has $\Delta^s = [\sigma]$,
thus we get a distinguished subdivision of $\overline{\Delta}$ by
setting $\overline{\Delta}^s = [\sigma]$ with the same centers of
subdivisions as in $\Delta$. Recall that to define the second
component of a distinguished pair, it suffices to construct the
embeddings of sheaves over the subfan $\overline{\Delta}^s =
[\sigma]$. For this we take the same embeddings as in $\Delta$.

Given the distinguished pair on $\overline{\Delta}$, we get the
Poincar\'e duality pairing  and the quadratic form $Q_{\overline{l}}$ on
$IH(\overline{\Delta})$. We define the evaluation maps on
$IH^{2n}(\Sigma)$ and $IH^{2n}(\overline{\Delta})$ using
the same metric on $\Lambda^n V^*$.

By the choice of the distinguished pairs on $\Sigma$ and
$\overline{\Delta}$, we have an identification $\cL_\Sigma(\Delta) =
\cL_{\overline{\Delta}}(\Delta)$.  Denote by $F$ the $A$-module 
\[ F = \{ (s_1, s_2) \in \cL(\Sigma)\times \cL(\overline{\Delta}) \mid
s_1|_\Delta = s_2|_\Delta\}.\]
On $F$ we define the degree $2$ map $l_F = (l, \overline{l})$, and the
quadratic form
\[ Q_F((s_1, s_2)) = Q_l (s_1) - Q_{\overline{l}} (s_2).\]
It is clear that $Q_F$ descends to a quadratic form on the quotient
$\overline{F} = F/ I F$.

\begin{lemma}\label{lem-def-beta} There exists a morphism of $A$-modules
\[ \beta: F \longrightarrow \cL(\hat{\Sigma}), \]
satisfying
\begin{enumerate}
\item[(a)] $\beta \circ l_F = \hat{l} \circ \beta;$
\item[(b)] $Q_F = Q_{\hat{l}}\circ \beta.$
\end{enumerate}
\end{lemma}

{\bf Proof.}  
Let $B = \Sym W^*$, where $V=W\oplus \RR v$ and $\Theta$ is a complete
fan in $W$. Using Lemma~\ref{lem-loc-prod} we have an
isomorphism of $A$-modules $\alpha:\cL(\Delta')\to \cL(\hat{\Delta})$
given as a composition:
\[ \cL(\Delta')\isom\cL(\Theta)\otimes_B A\isom\cL(\hat{\Delta}).\]
This isomorphism commutes with the actions of $\overline{l}$ and
$\hat{l}$ because both functions are defined as pullbacks of 
$l_\Theta$. By Lemma~\ref{lem-loc-prod1}, the isomorphism $\alpha$ is
also compatible with the representations of sections of $\cL$ as functions.

Define $\beta$ as follows. For $s = (s_1, s_2) \in F$,
set
\[ \beta(s)|_{\hat{\Sigma}-\hat{\Delta}} = s_1|_{\Sigma-\Delta},
\qquad \beta(s)|_{\hat{\Delta}} = \alpha (s_2|_{\Delta'}).\]
Since the restriction maps and $\alpha$ are $A$-linear, so is
$\beta$. The properties $(a)$ and $(b)$ follow easily from this
construction. We comment only on the minus sign in the expression of
$Q_F$. In Brion's construction of the evaluation map, if
$\sigma_1\in\hat{\Delta}$ and $\sigma_2\in\Delta'$ are maximal cones
that share a common facet, then the functions $\Phi_{\sigma_1}$ and
$\Phi_{\sigma_2}$ differ by a minus sign. Hence, $\alpha$ introduces a
minus sign into the evaluation map. \qed

We use the isomorphism $\alpha:\cL(\Delta')\to \cL(\hat{\Delta})$
defined in the proof of the previous lemma to construct a splitting
of the exact sequence
\[ 0\longrightarrow \cL(\Delta',\partial\Delta')\longrightarrow \cL(\overline{\Delta}) \longrightarrow
\cL(\Delta)\longrightarrow 0.\] 
If $s\in \cL(\Delta)$, we define $t\in \cL(\overline{\Delta})$ so 
that $t|_\Delta = s$ and $t|_{\Delta'} = \alpha^{-1} (s')$, where $s'$
is the image of $s$ under the embedding $\cL(\Delta)\subset
\cL(\hat{\Delta})$. Thus we get
\[ \cL(\overline{\Delta}) \isom \cL(\Delta) \oplus
\cL(\Delta',\partial\Delta').\]

\begin{lemma}\label{lem-ext-prim} Let $s_1\in \cL^{n-k}(\Sigma)$ be such that
  its image $[s_1] \in IH^{n-k}(\Sigma)$ is $l$-primitive. Then there
  exists an $s_2\in \cL^{n-k}(\overline{\Delta})$ satisfying:
\begin{enumerate}
\item[(a)] $(s_1,s_2) \in F$;
\item[(b)] the image $[s_2]\in IH^{n-k}(\overline{\Delta})$ is $\overline{l}$-primitive.
\end{enumerate}
\end{lemma}

{\bf Proof.} By Lemma~\ref{lem-deg2-isom} we have an isomorphism 
\[ IH(\Delta',\partial\Delta') \isom IH(\Theta)[-2],\]
commuting with the action of $\overline{l}|_{\Delta'}
= \pi^*(l_\Theta)$ and $l_\Theta$. 
By induction assumption, $l_\Theta$ defines a \LO on
$IH(\Theta)$ centered at degree $n-1$, hence
$\overline{l}$ defines a \LO on
$IH(\Delta',\partial\Delta')$ centered at degree $n+1$. 

Using the splitting 
\[ \cL(\overline{\Delta}) \isom \cL(\Delta) \oplus
\cL(\Delta',\partial\Delta'),\] 
we consider $s_1|_\Delta$ as an element in
$\cL(\overline{\Delta})$. Since the image of $l^{k+1} s_1$ in
$IH(\Sigma)$ is zero, we get that the image of $\overline{l}^{k+1}
s_1$ in $IH(\overline{\Delta})$ lies in the summand
$IH^{n+k+2}(\Delta',\partial\Delta')$. By the previous paragraph, 
\[ \overline{l}^{k+1}: IH^{n+1 -k-1}(\Delta',\partial\Delta') \longrightarrow 
IH^{n+1 +k+1}(\Delta',\partial\Delta') \]
is an isomorphism, hence there exists an element $s_0\in \cL^{n-k}
(\Delta',\partial\Delta')$, such that the image of $\overline{l}^{k+1}
(s_1|_\Delta+s_0)$ in $IH(\overline{\Delta})$ is zero. Thus $s_2 =
s_1|_\Delta+s_0$ is the required section. \qed

By the previous two lemmas we are reduced to proving that
$Q_{\overline{l}}$ satisfies the \HRM relations on
$IH(\overline{\Delta})$. Indeed,  given a primitive $[s_1] \in
IH^{n-k}(\Sigma)$, by Lemma~\ref{lem-ext-prim} there exists a primitive
$[s_2] \in IH^{n-k}(\overline{\Delta})$, such that $(s_1,s_2)\in F^{n-k}$. By
Lemma~\ref{lem-def-beta}(a) the class $[\beta(s_1,s_2)]\in
IH^{n-k}(\hat{\Delta})$ is $\hat{l}$-primitive, and by part $(b)$ of the same
lemma, we have
\[ (-1)^\frac{n-k}{2} Q_l([s_1]) = (-1)^\frac{n-k}{2} Q_{\hat{l}}
([\beta(s_1,s_2)]) +(-1)^\frac{n-k}{2} Q_{\overline{l}} ([s_2]).\]
The first summand on the right hand side is non-negative by induction
assumption, the second one is nonnegative if we know that
$Q_{\overline{l}}$ satisfies the \HRM relations, hence the left hand
side is nonnegative. It also
follows that the sum is zero if and only if both $[\beta(s_1,s_2)]$
and $[s_2]$ are zero;  this implies that in $\overline{F}$ we have
$[(s_1,s_2)] = [(s_1',0)]$ for some $s_1'$ such that the class
$[s_1']\in IH(\Sigma) \subset IH(\hat{\Sigma})$ is zero, hence
$[s_1]=[s_1']=0$. 

\begin{lemma}\label{lem-lef} The map $\overline{l}$ defines a \LO on
  $IH(\overline{\Delta})$ centered at degree $n$.
\end{lemma}

{\bf Proof.} By Poincar\'e duality it suffices to prove that
$\overline{l}^k:IH^{n-k}(\overline{\Delta})\to IH^{n+k}(\overline{\Delta})$ is
injective for $k>0$. If $h$ lies in the kernel of this map then by
Lemma~\ref{lem-rest} the restriction of $h$ to $IH(\Delta')$ is
zero. Considering the exact sequence
\[ 0\longrightarrow IH(\Delta,\partial\Delta)\longrightarrow IH(\overline{\Delta}) \longrightarrow
IH(\Delta')\longrightarrow 0,\] 
we have
that $h\in IH(\Delta,\partial\Delta)$. Since $\overline{l}$ acts on
$IH(\Delta,\partial\Delta)$ by multiplication with $l$, we see that
$h$ must lie in the kernel of $l^k$. We show in the next lemma that 
$l^k: IH^{n-k}(\Delta) \to IH^{n+k}(\Delta)$ is surjective, hence
dually $l^k:IH^{n-k}(\Delta,\partial\Delta)\to
IH^{n+k}(\Delta,\partial\Delta)$ is injective and $h=0$. \qed

\begin{lemma} The map 
\[ l^k: IH^{n-k}(\Delta) \longrightarrow IH^{n+k}(\Delta) \]
is surjective for $k> 0$.
\end{lemma}

{\bf Proof.} Let $p$ be the projection from the span of $\sigma$,
mapping $\Link_\Delta(\sigma)$ isomorphically to a complete fan
$\Lambda$ in $W\subset V$. Then we can write 
\[ \Delta = \Lambda \times_\phi [\sigma] \]
for a function $\phi:W\to \Span(\sigma)$. We also have the 
corresponding conewise linear isomorphism  
\[ \Phi: \Lambda \times [\sigma] \longrightarrow \Delta,\]
inducing an isomorphism of vector spaces
\[ \cL(\Delta) \isom \cL(\Lambda) \otimes \cL([\sigma]),\]
If we let $A_\Lambda$ be the ring of
polynomial functions on $W$, then the isomorphism of global
sections is an isomorphism of $A_\Lambda$ modules. 
Let $I_\sigma\subset \cA_\sigma$ be the ideal generated by linear
functions.

It suffices to show that $IH(\Delta)/l^k\cdot IH(\Delta)$ is zero in
degrees $n+k$ and higher. If we assume that $l|_\sigma=0$ then $l$ is
the pullback by $p$ of a strictly convex function on $\Lambda$. We have
\[ IH(\Delta)/l^k\cdot IH(\Delta) =  \cL(\Delta)/(I,l^k) = (IH(\Lambda)/(l^k)
\otimes \cL([\sigma]))/(I_\sigma),\]
where a linear function $f\in I_\sigma$ acts on the tensor product by
multiplication with $f\circ \Phi = 1\otimes f +f\circ\phi \otimes 1$.
Now applying the Hard Lefschetz theorem to $\Lambda$, we get
that $IH(\Lambda)/(l^k)$ is zero  in degrees $\dim\Lambda+k$ and
greater. Since $\cL([\sigma]) = \cL_\sigma$ has generators in 
degrees less than $\dim\sigma$, we get that $IH(\Delta)/l^k\cdot IH(\Delta)$
is zero in degrees $n = \dim\Lambda+k+\dim\sigma=n+k$ and higher. \qed

We divide the remaining proof that $Q_{\overline{l}}$ satisfies the
  \HRM relations into two cases.

\begin{lemma} Assume that $\dim \sigma = n$. Then there exists an
isomorphism 
\[ \gamma: IP_{\overline{l}} (\overline{\Delta}) \longrightarrow
IP_{l_\Theta}(\Theta),\]
such that $Q_{l_\Theta}(\gamma(p)) = Q_{\overline{l}}(p)$ for $p\in
IP_{\overline{l}} (\overline{\Delta})$. In particular,
$Q_{\overline{l}}$ satisfies the \HRM relations.
\end{lemma}

{\bf Proof.} Consider the splitting 
\[ IH(\overline{\Delta}) \isom IH(\Delta) \oplus
IH(\Delta',\partial\Delta').\] 
We have by construction of the primitive embeddings
\[ IH(\Delta) = IH([\sigma]) \isom IP_{l_\Theta}(\Theta).\]
We claim that $IH(\Delta)$ is the $\overline{l}$ primitive part of
$IH(\overline{\Delta})$, thus defining the isomorphism $\gamma$. To see
this, note that 
\[ IH(\Delta',\partial\Delta') \isom IH(\Theta)[-2], \]
and this isomorphism commutes with $\overline{l}$ and $l_\Theta$, thus
no element in $IH(\Delta',\partial\Delta')$ is primitive. On the other
hand, when we modify $\overline{l}$ so that it vanishes on $\sigma$,
it is clear that $\overline{l}$ maps $IH(\Delta)$ into
$IH(\Delta',\partial\Delta')$ by the composition:
\[ IH(\Delta) \isom IP_{l_\Theta}(\Theta) \stackrel{\iota[-2]}{\longrightarrow}
IH(\Theta)[-2] \isom IH(\Delta',\partial\Delta'),\]
where $\iota: IP_{l_\Theta}(\Theta) \to IH(\Theta)$ is the inclusion.
This shows that $IH(\Delta)= IP_{\overline{l}} (\overline{\Delta})$.

For $p\in IH^{n-k}(\Delta) \isom IP_{l_\Theta}(\Theta)$, we  represent
$p$ by a function $\pi^*(f)$. Then 
\[ Q_{\overline{l}}(p) = \langle \overline{l}^k \pi^*(f)\cdot
\pi^*(f)\rangle_{\overline{\Delta}} = \langle l_\Theta^{k-1} f\cdot
  f\rangle_\Theta  = Q_{l_\Theta} (\gamma(p)).\]
The middle equality follows from the fact that if $\overline{l}$
  vanishes on $\sigma$, then $\overline{l}(-v)>0$ and we can apply
  Lemma~\ref{lem-reduct}.
\qed

\begin{lemma}\label{lem-HRM} Assume that $\dim\sigma <n$. Then the
  quadratic form $Q_{\overline{l}}$  satisfies the \HRM 
  relations on $IH(\overline{\Delta})$.
\end{lemma}

{\bf Proof.} We write 
\[ \overline{\Delta} = \Lambda \times_\phi \overline{[\sigma]},\]
where $\overline{[\sigma]}\subset \overline{\Delta}$ is the subfan lying in the
linear span of $\sigma$,  $\Lambda$ is a fan conewise linearly
isomorphic to 
$\Link_{\overline{\Delta}}(\sigma)$, and $\phi$ is an appropriate
conewise linear map on $|\Lambda|$. Also let 
\[ \Phi:  \Lambda \times \overline{[\sigma]} \longrightarrow \overline{\Delta} \]
be the conewise linear isomorphism constructed from $\phi$. 

First we show that, changing $\overline{l}$ by a global linear
function if necessary, we may assume that $\Phi^* \overline{l}$ is
strictly convex on $\Lambda \times \overline{[\sigma]}$.
Let us assume that $\overline{l}$ vanishes on $\sigma$, hence
$\overline{l}|_\Delta$ is the pullback by the projection
$p:\Delta\to\Lambda$ of a strictly convex function $l_\Lambda$ on
$\Lambda$, and we write 
\[ \overline{l} = p^*{l_\Lambda} + (\overline{l}-p^*{l_\Lambda}),\]
where $(\overline{l}-p^*{l_\Lambda})$ is supported on $\Delta'$. Note
also that the restriction of $(\overline{l}-p^*{l_\Lambda})$ to
$\overline{[\sigma]}$ is strictly convex. Now 
\[ \Phi^* \overline{l} = \Phi^* p^*{l_\Lambda} + \Phi^*
(\overline{l}-p^*{l_\Lambda}),\] 
 hence $\Phi^* \overline{l}$ is the sum of pullbacks by the two
 projections of strictly convex functions on $\Lambda$ and on
 $\overline{[\sigma]}$. 

Let us consider a one-parameter family of skew products
\[ \overline{\Delta}_t = \Lambda \times_{\phi_t} \overline{[\sigma]},\]
where $\phi_t = t\phi$ for $0\leq t\leq 1$, and the corresponding
conewise linear isomorphisms
 \[ \Phi_t:  \Lambda \times \overline{[\sigma]} \longrightarrow \overline{\Delta}_t. \]
Then $\overline{\Delta}_1 = \overline{\Delta}$ and
$\overline{\Delta}_0 = \Lambda \times \overline{[\sigma]}$. We define the
strictly convex function $\overline{l}_t$ on $\overline{\Delta}_t$ by
\[ \overline{l}_t = (\Phi_t^{-1})^* \Phi^* (\overline{l}).\] 
Similarly, we pull back the distinguished pair on $\overline\Delta$ to
a distinguished pair on $\overline{\Delta}_t$ for all $0\leq t\leq 1$.

Note that the cohomology spaces $IH(\overline{\Delta}_t)$ have the
same dimension for all $t\in[0,1]$. By Lemma~\ref{lem-lef} applied to
the fan $\overline{\Delta}_t$, we also know that $\overline{l}_t$
induces  a \LO on $IH(\overline{\Delta}_t)$, hence the quadratic form
$Q_{\overline{l}_t}$ has maximal rank. Since all the data defining
$Q_{\overline{l}_t}$ is pulled back from the fan $\overline\Delta$, it
is clear that the signature of $Q_{\overline{l}_t}$ varies
continuously with t. Then 
since  $Q_{\overline{l}_0}$ has the correct signature by induction on
dimension and Corollary~\ref{cor-prod}, the same must
be true for any $t$.  \qed

\end{document}